\documentclass[graybox]{svmult}

\usepackage{type1cm}        
\usepackage{makeidx}         
\usepackage{graphicx}        
\usepackage{multicol}        
\usepackage[bottom]{footmisc}

\usepackage{booktabs}
\usepackage{multirow}
\usepackage{newtxtext}       %
\usepackage[varvw]{newtxmath}       

\makeindex             

\newcommand{\dealii}{\texttt{deal.II}}
\newcommand{\hp}{\textit{hp}}

\newcommand{\bu}{\mathbf{u}}
\newcommand{\bv}{\mathbf{v}}
\def\DIV{\nabla{\cdot\,}}
\def\SCAL{{\cdot}}
\def\GRAD{\nabla}

\title*{Solving the (Navier--)Stokes equations with space and time adaptivity using \dealii{}}
\author{Peter Munch\orcidID{0000-0003-2368-8533}, \\ Marc Fehling\orcidID{0000-0003-0984-793X}, \\ Martin Kronbichler\orcidID{0000-0001-8406-835X}, \\ Nils Margenberg\orcidID{0000-0003-2089-5545},  and\\ Laura Prieto Saavedra\orcidID{0000-0001-9086-2058}}
\authorrunning{P.\@ Munch et al.}
\institute{
Peter Munch \at Institute of Mathematics, Technical University of Berlin, Berlin, Germany \\\email{muench@math.tu-berlin.de}
\and Marc Fehling \at Faculty of Mathematics and Physics, Charles University, Prague, Czech Republic\\ \email{marc.fehling@matfyz.cuni.cz}
\and Martin Kronbichler \at Faculty of Mathematics, Ruhr University Bochum, Bochum, Germany\\ \email{martin.kronbichler@rub.de}
\and Nils Margenberg \at Faculty of Mathematics, Otto von Guericke University Magdeburg, Magdeburg, Germany \\\email{nils.margenberg@ovgu.de}
\and Laura Prieto Saavedra \at Weierstrass Institute for Applied Analysis and Stochastics, Berlin, Germany\\ \email{prieto@wias-berlin.de}}
\usepackage{subcaption}
\usepackage{bm}
\usepackage{tikz}
\usetikzlibrary{decorations}
\usetikzlibrary{arrows.meta}
\usetikzlibrary{decorations.pathreplacing, positioning, calc}
\usepackage{pgfplots}
\pgfplotsset{compat=1.18}

\graphicspath{{svg/},{svg_with_latex/},{svg_with_latex_script/}}

\begin{document}
\maketitle

\abstract*{TODO}

\abstract{In this article, we solve the Stokes and Navier--Stokes equations with the \dealii{} finite-element library. In particular, we use its multigrid, adaptive-mesh, and matrix-free infrastructures to design efficient linear and nonlinear iterative solvers, respectively. We solve the stationary Stokes equations on \hp-adaptive meshes with a \hp-multigrid approach, the transient Stokes equations with space-time finite elements and space-time multigrid, and, finally, the stabilized incompressible Navier--Stokes equations on locally refined meshes with a monolithic multigrid solver. The selected examples underline the flexibility and modularity of the multigrid infrastructure of \dealii{}.}

\section{Introduction}

Many solvers for finite element methods \textbf{(FEM)} are based on efficient solution methods for elliptic
partial differential equations \textbf{(PDE)}. Elliptic problems frequently occur in the context of, for example, computational fluid dynamics, geosciences, and computational solid mechanics. Efficient realizations often exploit the adaptivity of the mesh and of the local polynomial degree of the ansatz space, as well as efficient time-discretization schemes. The resulting linear or linearized system of equations needs efficient linear solvers, e.g., multigrid.

In this work, we summarize the key findings of three presentations held at ENUMATH 2025 in Heidelberg:  1) ``Efficient solvers for \hp-adaptive FEM computations'' by M.~Fehling, P.~Munch, M.~Kronbichler, W.~Bangerth; 2) ``An \hp-multigrid approach for tensor-product space-time finite element discretizations of the Stokes equations'' by N.~Margenberg, M.~Bause, P.~Munch; 3) ``Efficient distributed matrix-free multigrid methods for a stabilized solver for the incompressible Navier-Stokes equations on locally refined meshes'' by P.~Munch, L.~Prieto~Saavedra, B.~Blais. We solve three types of PDEs relevant in the context of computational fluid dynamics and geosciences:
\begin{itemize}
    \item stationary Stokes equations: find $\bu$ and $p$ such that
  \begin{align}
    \big( \nu \, \GRAD \bu,\ \GRAD \bv \big)
    - \big( p,\ \DIV \bv \big)
    +
    \big( \DIV \bu,\ q) & = 0 \quad \forall \bv,\, q;\label{eq:stokes}
  \end{align}
    \item transient Stokes equations: find $\bu$ and $p$ such that
  \begin{align}
    \left( \partial_t \bu ,\  \bv \right) 
    + \big( \nu \GRAD \bu,\ \GRAD \bv \big)
    - \big( p,\ \DIV \bv \big)
    + \big( \DIV \bu,\ q) & = 0 \quad \forall \bv,\, q;\label{eq:stokes_transient}
  \end{align}
    \item incompressible Navier--Stokes equations: find $\bu$ and $p$ such that
  \begin{align}
    F_{\text{NS}}=
    \left( \partial_t \bu ,\  \bv \right) 
    + \big( \bu \SCAL \GRAD \bu,\ \bv \big)
    + \big( \nu \GRAD \bu, \GRAD \bv \big)
    - \big( p,\ \DIV \bv \big)
    +
    \big( \DIV \bu,\ q) & = 0 \;\; \forall \bv,\, q.\label{eq:ns}
  \end{align}
\end{itemize}
We consider adaptivity in both space and polynomial degrees. For time discretization, we use BDF2 and tensor-product space-time FEM. The resulting linear or linearized systems of equations are solved with different variants of multigrid. For this purpose, we rely on the \dealii{} finite-element library~\cite{2025:arndt.bangerth.ea:deal, dealdesign} and especially on its multigrid infrastructure. All operators, particularly during smoothing and inter-grid transfer, are evaluated efficiently in a matrix-free approach.

The remainder of this work is organized as follows. In Section~\ref{sec:multigrid}, we give a short overview of the multigrid infrastructure in \dealii{}, focusing on its modularity and applicability to \hp-adaptivity. In Sections~\ref{sec:hp}--\ref{sec:stab}, we apply this infrastructure to different variants of the Stokes and Navier--Stokes equations. Finally, Section~\ref{sec:outlook} summarizes our conclusions and outlines directions for future research.

\section{Multigrid in \dealii{}}\label{sec:multigrid}

\dealii{} natively supports geometric, polynomial, and non-nested multigrid, and algebraic multigrid \textbf{(AMG)} is supported via wrappers to \texttt{PETSc}/\texttt{Trilinos}. Nesting these multigrid strategies gives hybrid multigrid algorithms. In particular, the multigrid infrastructure (and \dealii{} in general) allows users to work on locally refined meshes. Local refinement can occur both for the spatial mesh and in polynomial degree.

To support geometric, polynomial and non-nested multigrid in a unified framework, the implementation in \dealii{} must be modular and flexible. The crucial ingredient for a transfer operator is to implement an embedding between parent and child cells in the case of geometric multigrid and between cells and overlapping fine cells in the case of non-nested multigrid. In \dealii{}, this functionality is organized through a transfer base class, which is responsible for the transfer between two levels and can be composed (see \textbf{Fig.~\ref{fig:transfer}}). There are implementations of the transfer base class for geometric, polynomial, and non-nested multigrid, which work also for multivectors. The implementation of all of them follows the structure:
\begin{align}
\bm x^{(f)} = \mathcal{P}^{(f,c)} \circ \bm x^{(c)}
=
\sum_{e} \mathcal{S}^{(f)}_e \circ \mathcal{W}^{(f)}_e \circ \mathcal{P}^{(f,c)}_e \circ \mathcal{C}^{(c)}_e \circ \mathcal{G}^{(c)}_e \circ \bm x^{(c)},\label{eq:transfer}
\end{align}
where the prolongation is performed cell by cell. Here, we loop over all (coarse) cells, read the relevant degrees of freedom ($\mathcal{G}^{(c)}$) and resolve the constraints ($\mathcal{C}^{(c)}$), perform the prolongation locally ($\mathcal{P}^{(f,c)}$), and, finally, add contributions with some weighting to the (fine) result vector ($\mathcal{S}^{(f)}$, $\mathcal{W}^{(f)}$). In all algorithms in \dealii{}, the restriction operator is defined by the transpose (adjoint) of the prolongation operator. For more details, see  \cite{2023:munch.heister.ea:efficient}.

Just like in the case of transfer, smoothers and coarse-grid solvers can be flexibly chosen to best fit the application. For smoothing, \dealii{} provides, e.g., point Jacobi or additive Schwarz methods \textbf{(ASM)} as ingredients to relaxation smoothers or Chebyshev iterative schemes.
For the coarse-grid problem, both direct solvers and AMG are options provided by the software.

\newcommand{\R}{\textbf{R}}
\newcommand{\vct}[1]{\bm{#1}}
\newcommand{\symbf}[1]{\boldsymbol{#1}}
\newcommand{\symbfit}[1]{\boldsymbol{#1}}

\newcommand{\paran}[1]{\left(#1 \right)}
\newcommand{\norm}[1]{\left\lVert #1 \right\rVert}
\newcommand{\drv}{\mkern3mu\mathrm{d}}
\renewcommand{\coloneq}{\coloneqq}

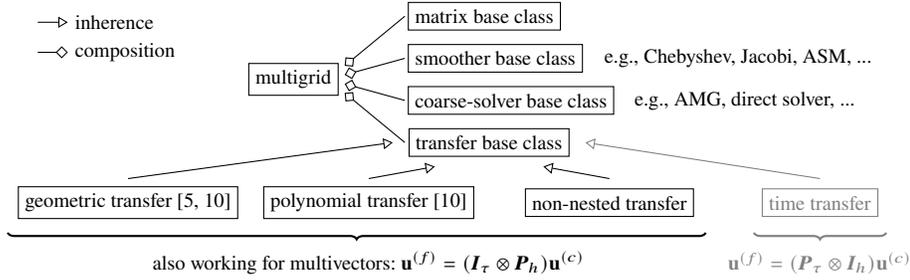
\begin{figure}[!t]

\centering

\begin{tikzpicture}[scale=0.8, every node/.style={scale=0.9}]
\node[draw] at (0,0) {transfer base class};

\node[draw] at (-6,-1.) {geometric transfer \cite{2021:clevenger.heister.ea:flexible, 2023:munch.heister.ea:efficient}};
\node[draw] at (-2,-1.) {polynomial transfer \cite{2023:munch.heister.ea:efficient}};
\node[draw] at (+2,-1.) {non-nested transfer};
\node[draw, gray] at (+5.5,-1.) {time transfer};


\draw [decorate, decoration = {brace,mirror}, thick] (-8,-1.5) --  (+3.6,-1.5);

\draw [decorate, decoration = {brace,mirror}, thick, gray] (+5.5-1.1,-1.5) --  (+5.5+1.1,-1.5);

\node at (-2,-2.0) {also working for multivectors: $\bu^{(f)}=
(\symbf{I}_\tau \otimes \symbf{P}_h) \bu^{(c)}$};

\node[gray] at (5.5,-2.0) {$\bu^{(f)}=
(\symbf{P}_\tau \otimes \symbf{I}_h) \bu^{(c)}$};

\draw[->, -{Triangle[open]}] (-6,-0.6) -- (-1.6,0);
\draw[->, -{Triangle[open]}] (-2,-0.6) -- (-0.9,-0.4);
\draw[->, -{Triangle[open]}] (+2,-0.6) -- (0.9,-0.4);
\draw[->, -{Triangle[open]}, gray] (+5.5,-0.6) -- (1.6,0);


\node[draw, anchor=west] at (-1.35,0.7) {coarse-solver base class};

\node[anchor=west] at (+2.3,0.7) {e.g., AMG, direct solver, ...};

\node[draw, anchor=west] at (-1.35,1.4) {smoother base class};

\node[anchor=west] at (+1.8,1.4) {e.g., Chebyshev, Jacobi, ASM, ...};

\node[draw, anchor=west] at (-1.35,2.1) {matrix base class};

\node[draw, anchor=east] at (-2.5,1.05) {multigrid};

\draw[->, -{Turned Square[open]}] (-1.4, 0.0) -- (-2.4,1.05-0.225);
\draw[->, -{Turned Square[open]}] (-1.4, 0.7) -- (-2.4,1.05-0.075);
\draw[->, -{Turned Square[open]}] (-1.4, 1.4) -- (-2.4,1.05+0.075);
\draw[->, -{Turned Square[open]}] (-1.4, 2.1) -- (-2.4,1.05+0.225);

\draw[->, -{Triangle[open]}] (-7-0.5,2) -- (-6.5-0.5,2);
\node[anchor=west] at (-6.5-0.5,2) {inherence};

\draw[->, -{Turned Square[open]}] (-7-0.5,1.5) -- (-6.5-0.5,1.5);
\node[anchor=west] at (-6.5-0.5,1.5) {composition};

\end{tikzpicture}

\caption{\dealii{} supports geometric, polynomial, and non-nested transfer, enabling multigrid variants with the same name. Furthermore, it supports multivectors. Time transfer is not supported, however, can be added in application codes, using inheritance.}\label{fig:transfer}

\end{figure}

\begin{figure}[!t]
    \centering
    \footnotesize\def\svgwidth{0.83\columnwidth}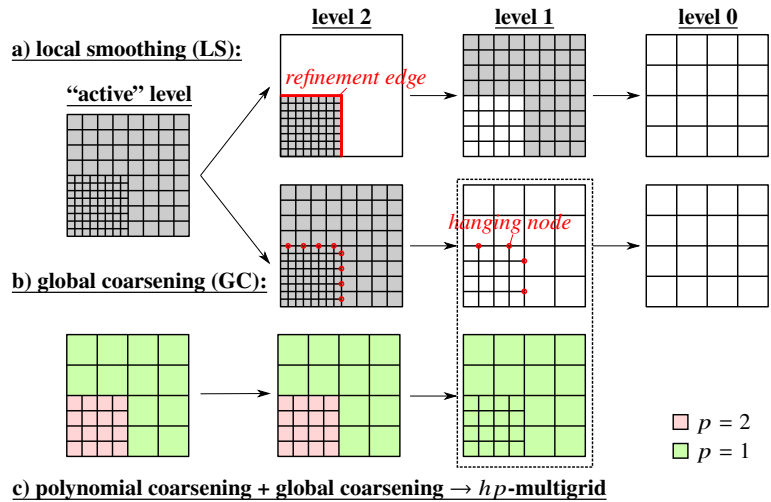

    \caption{Transfer types: (geometric) local smoothing, (geometric) global coarsening, and (global) polynomial coarsening. Nesting polynomial and geometric multigrid gives a (hybrid) \hp-multigrid.}\label{ref:vs}
\end{figure}

In the context of locally refined meshes, there are different strategies to create the levels:
\begin{itemize}
    \item \textbf{global coarsening (GC)}: smoothing is performed (globally) on the whole computational domain; the next level is obtained by coarsening (globally); this implies that one requires the resolution of (hanging-node) constraints during smoothing and transfer~\eqref{eq:transfer};
    \item \textbf{local smoothing (LS)}: smoothing is performed on the most refined part of the domain, resulting in internal domain boundaries, where homogeneous or inhomogeneous Dirichlet BCs are applied during pre- and postsmoothing; the next level is the next coarsest part of the domain.
\end{itemize}
These strategies are visualized for geometric multigrid in \textbf{Fig.~\ref{ref:vs}}a, b. In \dealii{}, $p$-multigrid and non-nested multigrid are always implemented in a global-coarsening way (see Fig~\ref{ref:vs}c).

Local smoothing and global coarsening for geometric multigrid were extensively compared in \cite{2023:munch.heister.ea:efficient} for the Poisson problem and a variable-viscosity Stokes simulation with a preconditioner based on block factorization. The main findings were that the number of iterations needed by LS and GC are similar, with a slight advantage for GC. In the serial setting, one iteration is slightly cheaper for LS, since fewer cells are processed and no hanging nodes have to be resolved. In contrast, GC tends to be more efficient when executed in parallel, since the work can be more straightforwardly balanced on each of the levels, resulting in less time spent waiting during synchronization. How well this is accomplished, we refer to as \textit{parallel workload efficiency} $\eta_w$. In contrast, \textit{vertical communication efficiency} $\eta_v$ quantifies how much data must be communicated during transfer. These performance metrics can be estimated from geometrical information and provide a first guideline for selecting an appropriate approach. The results in \cite{2023:munch.heister.ea:efficient} indicate that $\eta_w$ is the most significant metric.

\begin{figure}[!t]
\centering
\begin{subfigure}{0.32\textwidth}
\centering
\includegraphics[width=\textwidth]{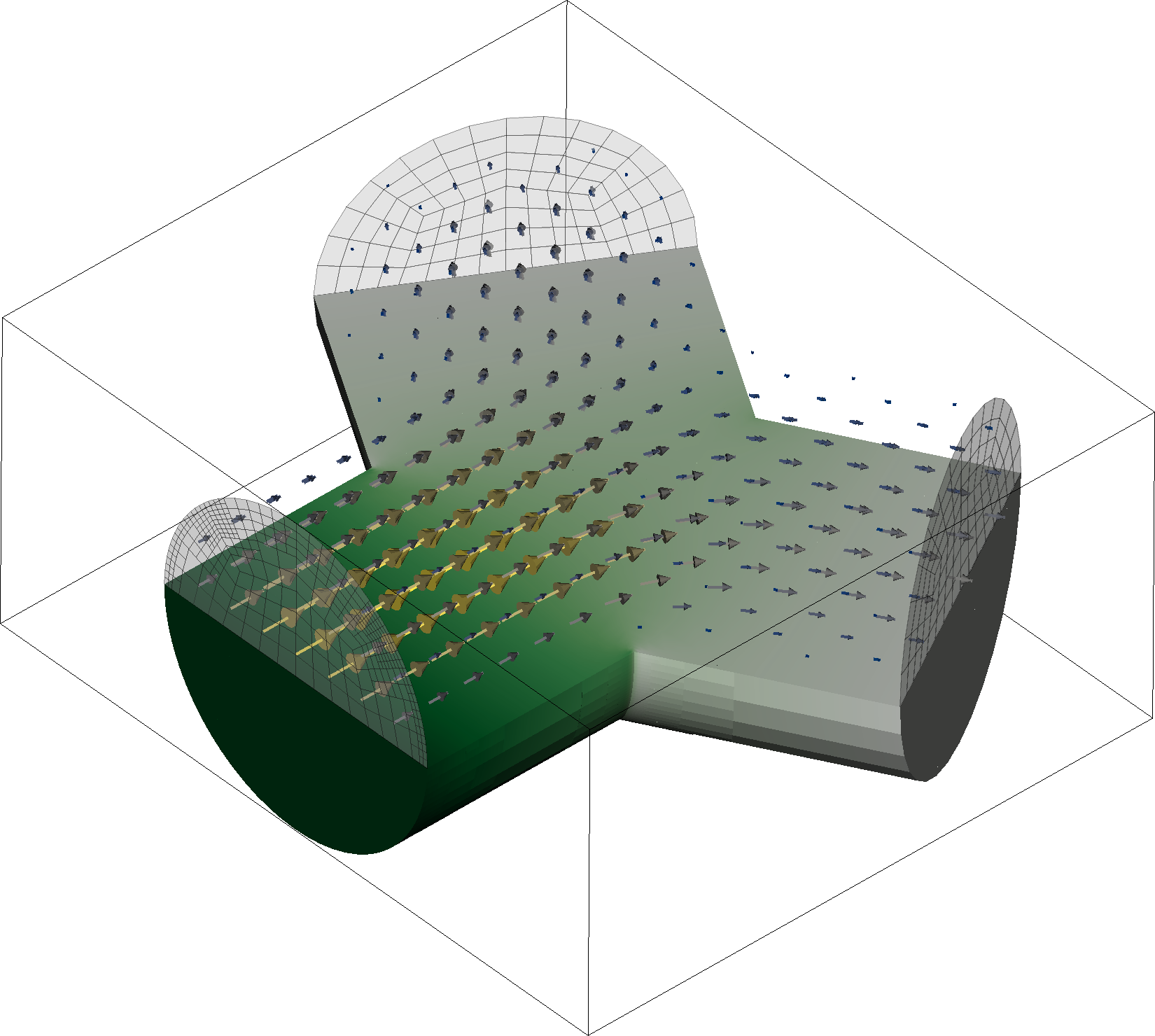}
\caption{Stokes flow through Y-pipe.}\label{fig:mesh_ypipe}
\end{subfigure}%
\hfill%
\begin{subfigure}{0.23\textwidth}
\centering
\includegraphics[width=\textwidth]{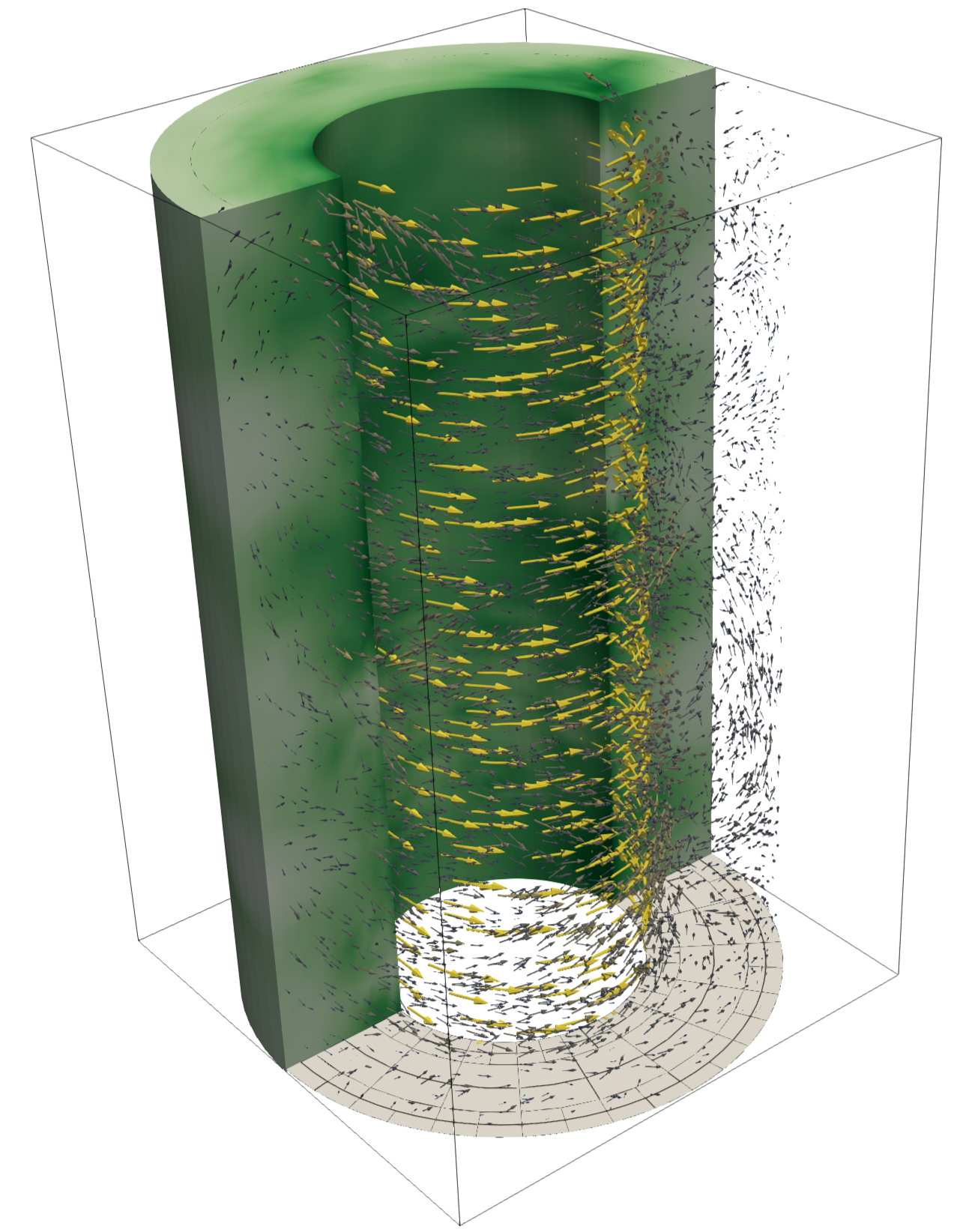}
\caption{Taylor--Couette flow.}\label{ref:mesh_taylor_couette}
\end{subfigure}%
\hfill%
\begin{subfigure}{0.42\textwidth}
\centering
\includegraphics[width=\textwidth]{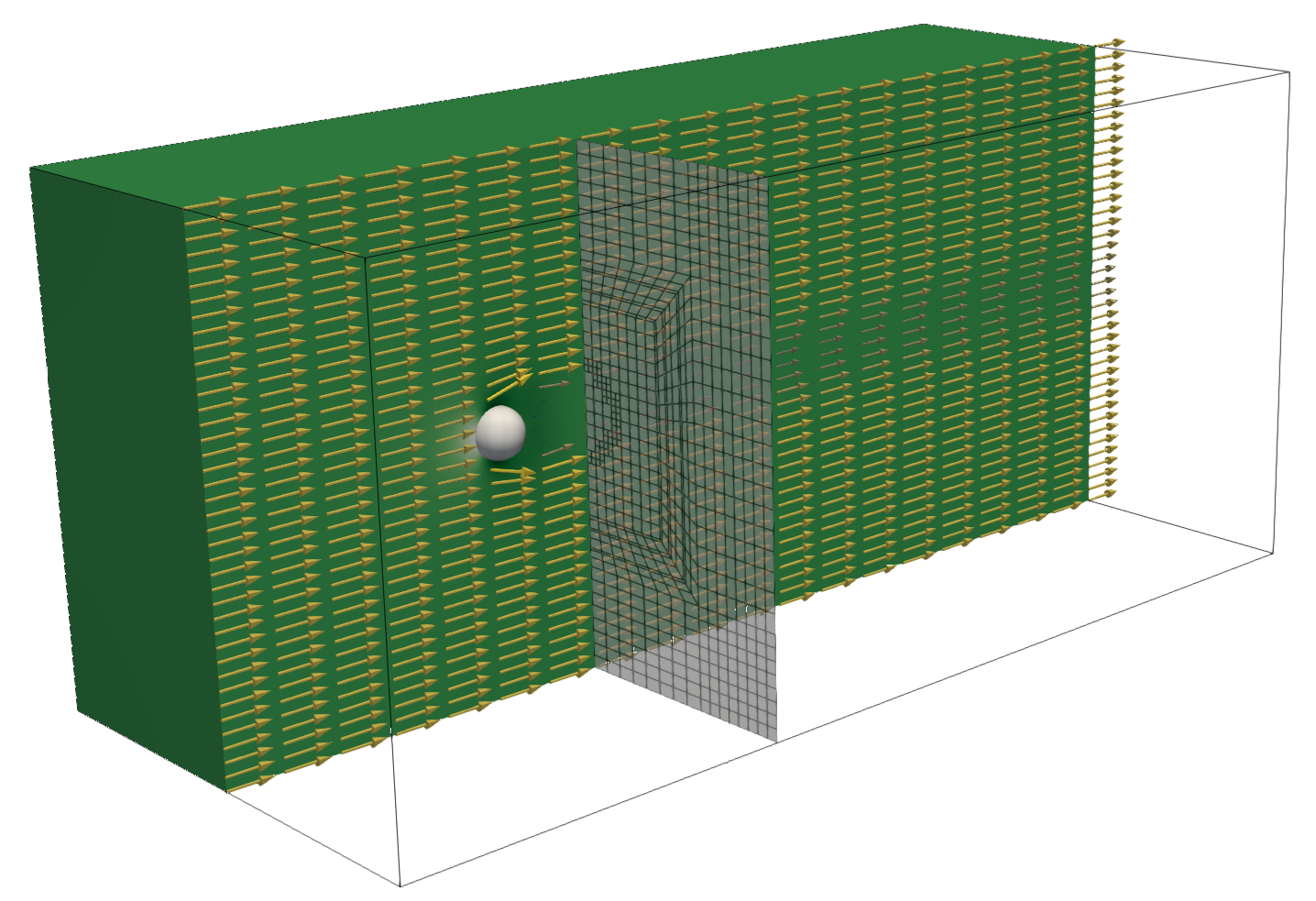}
\caption{Flow around a sphere.}\label{ref:mesh_sphere}
\end{subfigure}
\caption{Overview of the meshes used in experiments 1, 3, and 4. The color plots show the pressure and the vector plots the velocity.}\label{ref:geometries}
\end{figure}

\section{Multigrid for \hp-adaptive stationary Stokes}\label{sec:hp}

In this section, we solve \hp-adaptive stationary Stokes equations~\eqref{eq:stokes} with \hp-multigrid. Comprehensive support for parallel \hp-adaptive methods \cite{2023:fehling.bangerth:algorithms} has only recently been available in \dealii{}, and the task remains to find a proper preconditioner for linear systems arising from \hp-discretizations.
We present an elaborate approach using geometric multigrid for \hp-adapted problems.

For the Stokes system \eqref{eq:stokes}, we employ a block factorization and construct a block triangular, Silvester--Wathen type preconditioner 
\begin{align}
    A_{\text{Stokes}} &=
    \begin{bmatrix}
        A & B \\
        B^\top & 0 
    \end{bmatrix},
    &
    P_{\text{Stokes}} &=
    \begin{bmatrix}
        A & B \\
        0 & -S 
    \end{bmatrix},
    &
    P_{\text{Stokes}}^{-1} &=
    \begin{bmatrix}
        A^{-1} & A^{-1} B S^{-1} \\
        0 & - S^{-1}
    \end{bmatrix},
    \label{eq:stokes_block}
\end{align}
with the viscous vector Laplacian $A$, the discrete gradient $B$, and the negative Schur complement $S=B^\top A^{-1} B$. We approximate $S$ by a scaled mass matrix and compute its inverse with the conjugate gradient method preconditioned by point Jacobi. We solve the Stokes system with FGMRES. Furthermore, we approximate the action of $A^{-1}$ by a single V-cycle of \hp-multigrid.

\hp-multigrid has already been realized for discontinuous elements \cite{2026:nastase:high-order} and continuous elements with $h$-hierarchical bases \cite{2010:mitchell:hp-multigrid}, but herein we show a more general approach to support a larger set of continuous elements. In \cite{2021:pazner:uniform}, an efficient two-level preconditioner based on a low-order refined mesh was proposed for the DG case.

In essence, our algorithm is a generalization of \cite{2023:munch.heister.ea:efficient} for \hp-adapted discretizations. In the multigrid hierarchy, first the polynomial degree $p$ of the elements is decreased one by one until a discretization with all linear elements is obtained. Next, (geometric) global coarsening is applied until the coarse grid is reached (see Fig.~\ref{ref:vs}b+c). We resolve hanging-node constraints for \hp-discretizations in the case of continuous elements \cite{2009:bangerth:data}. 

As smoother, we use a Chebyshev iteration with degree 5 around either point Jacobi or a blend of point Jacobi and ASM. The blocks of ASM comprise unknowns on faces where one neighbor cell is $h$-refined compared to the neighbor and the neighbor is $p$-refined compared to the cell. The standard process of eliminating the constrained degrees of freedom~\cite{2009:bangerth:data} in this configuration leads to ill-conditioned relations between the remaining unknowns, translating to a bad conditioning of the linear system to be solved, as indicated by the eigenvalue estimates in \textbf{Fig.~\ref{fig:hp_comparison_iterations}}. For unknowns that are not part of any ASM blocks, the ASM strategy employs point Jacobi. As coarse-grid solver, we use AMG.


\textbf{Exp.~1}: We solve a 3D Stokes flow problem through a Y-pipe (\textbf{Fig.~\ref{fig:mesh_ypipe}}, \cite{2023:fehling.bangerth:algorithms}). As boundary conditions, we impose Hagen--Poiseuille inflow on one opening, zero-traction outflow on the others, as well as no-slip conditions on the lateral surfaces. We perform multiple successive solve-and-refine cycles, during each  10\,\% of cells with the largest error determined by the Kelly error estimator are refined, of which in turn half are $p$-refined based on estimated smoothness by the decay of Legendre coefficients. Both estimators are readily available in \dealii{}. We use a selection of stable Taylor--Hood continuous elements $\bm{Q}_p/Q_{p-1}$ with $p \in [3,8]$. 
We vary the multigrid type and compare the performance between a matrix-based and matrix-free implementation. For the former, we use AMG to approximate the action of $A^{-1}$.

\begin{figure}[!t]
\centering
\begin{tikzpicture}
\small
\begin{axis}[
  scale=0.725,
  xlabel={Refinement cycle},
  ylabel={Number of outer solver iterations},
  xlabel near ticks,
  ylabel near ticks,
  xtick distance=1,
  legend style={at={(0.97, 0.35)}, anchor=east},
  legend cell align=left,
]
\addplot+ table [y=iterations, x=cycle, col sep=comma] {data/ypipe/stokes_amg_cycles_256.csv};
\addlegendentry{AMG};
\addplot+ table [y=iterations, x=cycle, col sep=comma] {data/ypipe/stokes_mg-diagonal_cycles_256.csv};
\addlegendentry{\hp-MG: Jacobi};
\addplot+ table [y=iterations, x=cycle, col sep=comma] {data/ypipe/stokes_mg-asm_cycles_256.csv};
\addlegendentry{\hp-MG: Jacobi+ASM};
\end{axis}
\end{tikzpicture}
\hfill%
\begin{tikzpicture}
\small
\begin{axis}[
  scale=0.725,
  xlabel={Refinement cycle},
  ylabel={Maximum eigenvalue of $P^{-1}A$},
  xlabel near ticks,
  ylabel near ticks,
  xtick distance=1,
  legend pos=north west,
  legend cell align=left,
]
\pgfplotsset{cycle list shift=1};
\addplot+ table [y=max_ev, x=cycle, col sep=comma] {data/ypipe/stokes_mg-diagonal_cycles_256.csv};
\addlegendentry{$P^{-1}$: Jacobi};
\addplot+ table [y=max_ev, x=cycle, col sep=comma] {data/ypipe/stokes_mg-asm_cycles_256.csv};
\addlegendentry{$P^{-1}$: Jacobi+ASM};
\end{axis}
\end{tikzpicture}
\caption{Flow through Y-pipe (Exp. 1): robustness metrics after solving the Stokes system with successive \hp-adaptive refinements on 256 MPI processes. Left: comparison of the number of outer solver iterations, using AMG or \hp-multigrid with plane or ASM-enhanced point Jacobi for smoothing. Right: estimation of maximum eigenvalue of the preconditioned $A$-block.}
\label{fig:hp_comparison_iterations}
%
\vspace{1em}
\centering
\begin{tikzpicture}
\small
\begin{loglogaxis}[
  scale=0.725,
  xlabel={Number of degrees of freedom (DoFs)},
  ylabel={Wall-clock time [seconds]},
  xlabel near ticks,
  ylabel near ticks,
  legend cell align=left,
  legend pos=outer north east
]
\addplot+ table [y=solve_max, x=dofs, col sep=comma] {data/ypipe/stokes_amg_cycles_256.csv};
\addlegendentry{AMG};
\addplot+ table [y=solve_max, x=dofs, col sep=comma] {data/ypipe/stokes_mg-diagonal_cycles_256.csv};
\addlegendentry{\hp-MG w.\@ Jacobi};
\addplot+ table [y=solve_max, x=dofs, col sep=comma] {data/ypipe/stokes_mg-asm_cycles_256.csv};
\addlegendentry{\hp-MG w.\@ Jacobi+ASM};
\addplot[thick, samples=2, domain=3365108:117269100] {5e-7*x};
\addlegendentry{$\mathcal{O}(N_\text{DoFs})$};
\addplot[thick, dashed, samples=2, domain=3365108:27470962] {5e-13*x*x};
\addlegendentry{$\mathcal{O}^2(N_\text{DoFs})$};
\end{loglogaxis}
\end{tikzpicture}
\caption{Flow through Y-pipe (Exp. 1): comparison of wall-clock times on 256 MPI processes. We perform successive \hp-adaptive refinements and solve the Stokes system with different types of preconditioners applied to the $A$-block.}
\label{fig:hp_comparison_time}
\end{figure}

We solve the Stokes system on the Expanse supercomputer \cite{2021:strande:expanse}. Each computing node consists of two AMD EPYC\texttrademark{} 7742 processors (64 cores each, 2.25\,GHz) and 256\,GB of DDR4 memory (16 memory channels). Nodes are connected with a Mellanox\textsuperscript{\textregistered} HDR-100 InfiniBand Interconnect network operating in a hybrid fat-tree topology. Experiments were performed on 2 nodes with a total of 256 MPI processes.

A comparison of solver iterations in \textbf{Fig.~\ref{fig:hp_comparison_iterations}} shows that \hp-multigrid is more robust compared to AMG. With the ASM-enhanced smoother, we observe a number of $28\pm2$ iterations independent of the level of refinement, which is the same number for pure $h$-refinement with fixed $p$. In contrast, the number of iterations for point Jacobi increases with the number of refinement cycles. In \textbf{Fig.~\ref{fig:hp_comparison_time}}, we found the runtime of AMG to behave roughly as $\mathcal{O}(N^2)$ and the runtime of \hp-multigrid as $\mathcal{O}(N)$, where $N$ is the number of degrees of freedom. The ASM-enhanced smoother method has a slight advantage compared to point Jacobi, especially for the higher degrees, because the increased cost per iteration is offset by a robust behavior of iteration counts.

Overall, we see that solvers with \hp-adaptivity, using \hp-multigrid, are indeed comparable regarding robustness and scalability to equivalent solvers for pure $h$-adaptive discretizations. We will investigate their suitability for transient problems in the future.

\section{Space-time multigrid for transient Stokes}\label{sec:stfem}

In this section, we solve the transient Stokes equations~\eqref{eq:stokes_transient} with space-time multigrid.
On each time slab $I_n=(t_{n-1},\,t_n]$ with $\tau_n=t_n-t_{n-1}$, we use a $Q$-stage variational time discretization
(e.g.\ DG($k$) with $Q=k+1$, or collocation implicit Runge--Kutta (IRK) with tableau $(\mathbf G_Q,\,\mathbf b_Q,\,\mathbf c_Q)$) and work with stage-wise Stokes vectors
\[
  \mathbf x_{n,q}:=(\mathbf u_{n,q},\, \mathbf p_{n,q})^{\top}\in \mathbb{R}^{N_u+N_p},
\qquad q=1,\dots,Q,
\]
collected in the \emph{stage-major} block vector
\(
  \mathbf x_n := (\mathbf x_{n,1},\dots,\,\mathbf x_{n,Q})^{\top}\in \mathbb{R}^{Q(N_u+N_p)}.
\)
The variational time discretization (see~\cite{2025:margenberg.bause.ea:hp}) yields, on $I_n$, the raw Kronecker system
\(
\mathbf K_\square^{\mathrm{raw}} = \mathbf A_Q^\tau\otimes \mathbf M_\square + \tau_n\,\mathbf M_Q\otimes \mathbf A_\square,
\)
with purely temporal matrices $\mathbf M_Q$ (mass) and $\mathbf A_Q^\tau$ (stiffness, incl.\ the DG jump term), see
\cite{2024:margenberg.munch:space-time, 2025:margenberg.bause.ea:hp}. Hence, temporal operators act as dense $Q\times Q$ matrices on the
stage index, while spatial blocks act stage-wise. Mass scaling (left multiplication by $M_Q^{-1}$) defines
$\mathbf G_Q:=\mathbf M_Q^{-1}\mathbf A_Q^\tau$ and gives, in stage-major ordering,
\begin{equation}
  \label{eq:proceedings-kronecker}
  \mathbf K_\square \mathbf x_n = (\mathbf G_Q\otimes \mathbf M_\square + \tau_n\,\mathbf I_Q\otimes \mathbf A_\square)\mathbf x_n=\mathbf b_n,
  \quad
  \mathbf b_n=\mathbf f_n+(\mathbf m_Q\otimes \mathbf M_{\square})x_{n-1,q},
\end{equation}
where $\mathbf f_n$ is the stage load vector and $m_Q$ encodes the jump contribution
from the previous slab (see \cite{2024:margenberg.munch:space-time,
  2025:margenberg.bause.ea:hp} for $\mathbf M_Q,\,\mathbf A_Q^\tau,\,\mathbf m_Q$). The
operator~\eqref{eq:proceedings-kronecker} has the same IRK-type stage Kronecker
structure: in collocation IRK, the stage matrix is determined by the Butcher
tableau, whereas, for DG($k$), it is induced by the temporal basis, quadrature,
and jump term. For the instationary Stokes equations, we take
\(\mathbf M_\square\) as the Stokes block mass operator (velocity mass, zero pressure block) and set
\(\mathbf A_\square \equiv \mathbf A_{\mathrm{Stokes}}\) (see~\eqref{eq:stokes_block}). We use inf--sup stable
elements in space, e.\,g., $\bm{Q}_{p}^d/P_{p-1}^{\mathrm{disc}}$ (and, in principle, other pairs supported by \dealii). The resulting
monolithic stage system couples velocity and pressure across all temporal stages and is generally non-symmetric due to the temporal
derivative contribution through $\mathbf A_Q^\tau$.

We precondition GMRES by a single V-cycle of an \hp{} space-time multigrid method (\hp-STMG), applied to the monolithic
operator~\eqref{eq:proceedings-kronecker}. The multigrid hierarchy is built on a tensor-product
space-time FE space and combines $p$- and $h$-coarsening
in both space and time. Two rules determine the order in which the
levels are ordered:
\textbf{R1}:~\emph{Coarsening in space before time} (geometric coarsening first in $h$, then in $\tau$); \textbf{R2}:~\emph{Coarsening polynomially before geometrically} (halve $p$/$k$ before coarsening $h$/$\tau$).

In a time-marching setting, we omit geometric coarsening in time and retain only $p$-coarsening
in time (while still using \hp-coarsening in space).
Let $\boldsymbol P_h,\boldsymbol R_h$ denote the spatial prolongation/restriction and $\boldsymbol P_\tau,\boldsymbol R_\tau$ their temporal analogues. The natural
space-time transfers are tensor products
\[
\boldsymbol P = \boldsymbol P_\tau\otimes \boldsymbol P_h,
\qquad
\boldsymbol R = \boldsymbol R_\tau\otimes \boldsymbol R_h,
\]
which allows us to apply them sequentially in the implementation:
\[
\bu^{(f)}=(\boldsymbol I_Q\otimes \boldsymbol P_h)\,\bu^{(c)},\qquad
\bu^{(f)}=(\boldsymbol P_\tau\otimes \boldsymbol I_h)\,\bu^{(c)},
\]
and analogously for restriction. This preserves modularity and let us
use the \dealii{} infrastructure. The spatial transfers $\boldsymbol P_h,\, \boldsymbol R_h$ are taken
directly from \dealii{}, while temporal transfers are not provided by the
library. We implement $\boldsymbol P_\tau,\, \boldsymbol R_\tau$ as custom transfer operators that inherit
from the \dealii{} transfer base class (Fig.~\ref{fig:transfer}) and can be
plugged into the multigrid driver like native transfers. Concretely, $\boldsymbol P_\tau$
and $\boldsymbol R_\tau$ are implemented as linear combinations of the stage values at
a spatial degree of freedom.

We use a space-time ASM smoother that resolves the
velocity-pressure coupling locally on patches. On each patch (one spatial cell over one
time slab), we gather the associated space-time velocity and pressure DoFs. The implementation also supports vertex-star
patches, while keeping the temporal direction cell-based.

On the coarsest level (small $p,k$ and coarse $h,\tau$), we solve directly. The
\hp-STMG hierarchy and the patch smoother are implemented in a matrix-free fashion. At high
order, the cost is dominated by the ASM smoother, motivating future work.

\textbf{Exp.~2:} We consider a variation of the benchmark ``flow around a cylinder'' by
Sch\"afer and Turek~\cite{SchaeferTurek1996}, restricted to Stokes. The flow is driven by a Dirichlet
inflow profile $\bu=\bu^D$ on $\Gamma_{\mathrm{in}}= 0\times[0,\,H]$ (channel height $H=0.41$),
\begin{equation}
\bu^D(x,y,t) = \sin\left(\frac{\pi}{8} t\right)\frac{6y(H-y)}{H^2}\,.
\end{equation}
We impose $v=0$ on $\Gamma_{\mathrm{wall}}$
and the obstacle, and a do-nothing outflow condition on
$\Gamma_{\mathrm{out}}$. In \textbf{Table~\ref{tab:gmres-k-h}}, we show iteration numbers and wall-clock times under mesh refinement for $Q_{k}^d/P_{k-1}^{\mathrm{disc}}$ elements with a DG($k-1$) time discretization.
All experiments were run on a single node (HSUper at Helmut Schmidt University) using 72 MPI ranks.
The node is equipped with $2\times$ Intel Xeon Platinum 8360Y CPUs (36 cores each) and 256\,GB RAM.

\begin{table}[t]
\centering
\caption{Flow past cylinder (Exp.~2): Average number of GMRES iterations $\overline{N}_L$ and total wall-clock time $T$ under uniform $h$-refinement for increasing polynomial degree $k$.}
\label{tab:gmres-k-h}
\begin{tabular}{c r r r@{\hskip 5pt}c r@{\hskip 5pt}c r@{\hskip 5pt}c r@{\hskip 5pt}c}
\toprule
\multirow{2}{*}{ref.} & \multirow{2}{*}{cells} & \multirow{2}{*}{$N_T$}
& \multicolumn{2}{c}{$k=1$}
& \multicolumn{2}{c}{$k=2$}
& \multicolumn{2}{c}{$k=3$}
& \multicolumn{2}{c}{$k=4$} \\
\cmidrule(lr){4-5}\cmidrule(lr){6-7}\cmidrule(lr){8-9}\cmidrule(lr){10-11}
& & & $\overline{N}_L$ & $T$ [s]
    & $\overline{N}_L$ & $T$ [s]
    & $\overline{N}_L$ & $T$ [s]
    & $\overline{N}_L$ & $T$ [s] \\
\midrule
1 &  1600 &  849  & 11.74 & $2.33$                 & 17.12 & $2.11\!\cdot\!10^{1}$ & 16.69 & $5.99\!\cdot\!10^{1}$ & 13.56 & $1.02\!\cdot\!10^{2}$ \\
2 &  6400 & 1697  &  8.77 & $1.25\!\cdot\!10^{1}$  & 10.16 & $9.81\!\cdot\!10^{1}$ & 11.32 & $3.14\!\cdot\!10^{2}$ & 11.37 & $6.82\!\cdot\!10^{2}$ \\
3 & 25600 & 3393  &  7.57 & $8.96\!\cdot\!10^{1}$  &  7.64 & $5.97\!\cdot\!10^{2}$ &  9.28 & $2.06\!\cdot\!10^{3}$ &  9.88 & $4.74\!\cdot\!10^{3}$ \\
\bottomrule
\end{tabular}
\end{table}
The GMRES wall-clock time follows the standard cost model
$T \approx N_T\,  \overline{N}_L\, \overline{T}_{L}$ and, therefore, grows mainly with the number of time steps $N_T$ and the per-iteration cost $\overline{T}_{L}$, which increases under $h$-refinement (larger global dimension) and with polynomial degree (more DoFs per cell).
Accordingly, from ref.~1 to ref.~3, $T$ increases by $\approx 38\times$ for $k=1$ and by $\approx 46\times$ for $k=4$. For fixed refinement, raising $k$ from 1 to 4 yields $\approx 44\times$ (ref.~1) up to $\approx 53\times$ (ref.~3).
In contrast, the average number of linear iterations $\overline{N}_L$ stays $\mathcal{O}(10)$, suggesting robust preconditioned convergence: it decreases with refinement for $k\le 2$ and is essentially mesh-independent for $k\ge 3$ (about $9$--$11$).

Future work will focus on reducing the cost of the ASM smoother,
in particular by developing efficient multilevel block preconditioners and
cheaper patch solvers.

\section{Multigrid for stabilized Navier--Stokes on locally refined meshes}\label{sec:stab}

In this section, we solve the incompressible Navier--Stokes equations considering two types of stabilization: Streamline-Upwind Petrov-Galerkin (SUPG) and Pressure-Stabilized Petrov-Galerkin (PSPG). These techniques require two terms to be added to the Navier--Stokes equations~\eqref{eq:ns} as follows:
 find $\bu$ and $p$ such that
\begin{align*}
    F_{\text{NS}} + F_{\text{SUPG}} + F_{\text{PSPG}} = 0 \quad \forall \bv,\, q,
\end{align*}
with
\begin{align*}
 F_\text{SUPG} &= \sum_k \delta_1 \big( \partial_t \bu + \bu \SCAL \GRAD \bu + \GRAD p - \nu \Delta \bu,\ \bu \cdot \GRAD \bv \big)_{\Omega_k},   \\
 F_\text{PSPG} &= \sum_k \delta_1 \big( \partial_t \bu + \bu \SCAL \GRAD \bu + \GRAD p - \nu \Delta \bu,\ \GRAD q \big)_{\Omega_k}.
\end{align*}
The SUPG terms reduce oscillations in the velocity field when simulating flows with high Reynolds numbers, while PSPG allows using equal-order elements for the velocity and the pressure $\bm{Q}_p Q_p$. In time, we discretize the system, using second-order backward differentiation (BDF2). The resulting system of nonlinear equations is solved with a Newton–Krylov method, where the Jacobian is approximated by one V-cycle of a monolithic hybrid multigrid (first: geometric coarsening; then: polynomial coarsening). The side effect of the PSPG stabilization is that the pressure diagonal block is non-zero in the Jacobian $J$ of the system:
\begin{align*}
    J=\begin{bmatrix}
        A' & B' \\
        C & D
    \end{bmatrix},
\end{align*}
allowing us to use simpler smoothers, e.g., point Jacobi. The solver is part of the \texttt{Lethe} library and its efficiency for globally refined meshes was shown in \cite{2025:prieto-saavedra.munch.ea:matrix-free}. Here, we consider two three-dimensional cases with locally refined meshes:
\begin{itemize}
    \item \textbf{Exp. 3}: Taylor--Couette flow at $\mathrm{Re}=4000$; transient flow between coaxial cylinders; static outer cylinder and inner cylinder with a fixed angular velocity; static local refinement of cells close to the curved wall (see Fig.~\ref{ref:mesh_taylor_couette}); fixed \texttt{CFL} equal to 1; multigrid using relaxation method with 5 iterations around point Jacobi as a smoother and a direct coarse-grid solver.
    \item \textbf{Exp. 4}: stationary flow around a sphere at $\mathrm{Re}=150$ with a fixed entrance velocity; region of mesh refinement is determined by a Kelly error estimator based on the pressure (see Fig.~\ref{ref:mesh_sphere}); initial condition obtained by ramping up the $\mathrm{Re}$ number starting with $\mathrm{Re}=10$; multigrid using relaxation method with 2 iterations around ASM as a smoother and a direct coarse-grid solver.
\end{itemize}
We compare both local-smoothing and global-coarsening (geometric) multigrid algorithms. We perform experiments of the two cases, using linear/quadratic elements ($\bm{Q}_1Q_1$/$\bm{Q}_2Q_2$) and different initial global refinements~$l$. Adopting the local-smoothing algorithm without modifications implies that we have Dirichlet BCs both for $\bu$ and $p$ at refinement edges. The experiments are conducted on 12 cores of an Intel(R) Core(TM) i9-14900 machine. \textbf{Tables~\ref{tab:tc}} and \textbf{\ref{tab:fpc}} show the performance metrics, relevant iteration numbers, and simulation times. 

\begin{figure}[!t]
    \centering

    \includegraphics[width=0.9\textwidth]{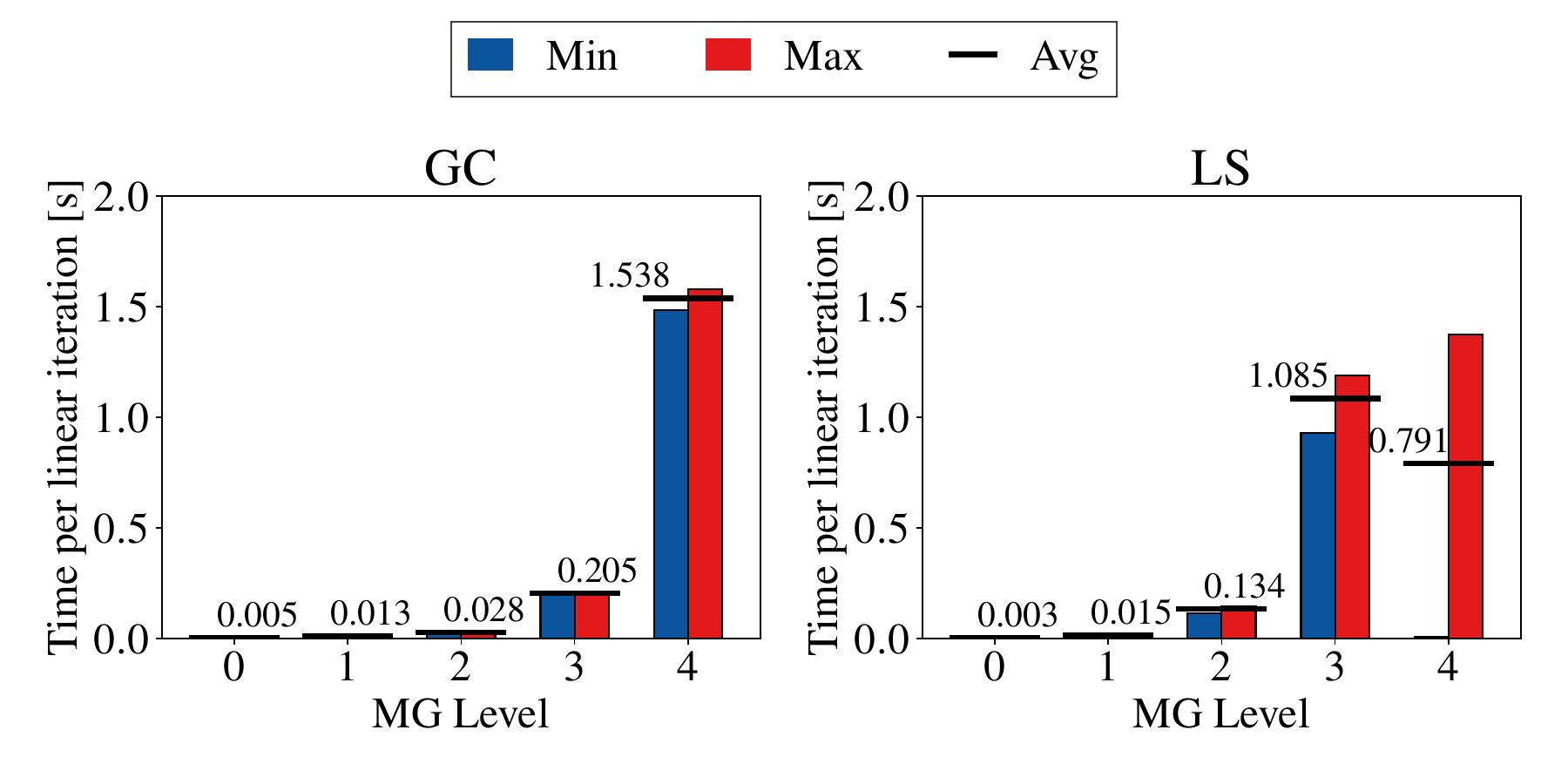}

    \caption{Flow around a sphere (Exp.~4): time per linear iteration on each level for the $\bm{Q}_2Q_2$ case with $l=2$. It showcases the average, minimum, and maximum time among all processes.}\label{fig:fpc:levels}
\end{figure}

The results (see Tables~\ref{tab:tc} and \ref{tab:fpc}) show that the simulation times for GC and LS are similar for Exp.~3, with a slight advantage of GC over LS; this advantage is more pronounced in Exp.~4, where the times are up to 30\% lower for GC. The reason is that the efficiencies $\eta_w$, $\eta_v$ are comparable for Exp.~3, while, in Exp.~4, $\eta_w$ is much worse for LS, indicating load imbalances between the cores on different levels. These load imbalances increase the time of a single iteration, as indicated in \textbf{Fig.~\ref{fig:fpc:levels}}, where, for $\bm{Q}_2Q_2$, the min./max./avg. time per linear iteration on each level spent among all processes is shown. The results are consistent with the ones obtained for a simple Poisson operator in \cite{2023:munch.heister.ea:efficient}.

\begin{table}[!t]
    \centering

    \caption{Taylor--Couette flow (Exp. 3): number of transient iterations $N_T$, number of nonlinear iterations per transient iteration $\overline{N}_N$, number of linear iterations per nonlinear iteration $\overline{N}_L$, and wall-clock time $T$. Performance metrics for the cases with $l=3$: $\eta_w$ parallel workload efficiency and $\eta_v$ vertical communication efficiency.}\label{tab:tc}

\begin{tabular}{c|cc|cc|cc}
\hline
 & \multicolumn{6}{c}{$Q_1Q_1$} \\ \hline
 & \multicolumn{2}{c|}{} & \multicolumn{2}{c|}{GC} & \multicolumn{2}{c}{LS} \\ \hline
$l$ & $N_T$ & $\overline{N}_N$ & $\overline{N}_L$ & $T$  & $\overline{N}_L$ & $T$ \\ \hline
2 & 72 & 2 & 4.71 & 6.47  & 7.02 & 7.97\\
3 & 152 & 2 & 4.01   & 84.5   & 5.03 & 84.5\\\hline
\end{tabular}\quad
\begin{tabular}{c|cc|cc|cc}
\hline
 & \multicolumn{6}{c}{$Q_2Q_2$} \\ \hline
 & \multicolumn{2}{c|}{} & \multicolumn{2}{c|}{GC} & \multicolumn{2}{c}{LS} \\ \hline
$l$ & $N_T$ & $\overline{N}_N$ & $\overline{N}_L$ & $T$ & $\overline{N}_L$ & $T$ \\ \hline
2 & 149 & 2    & 4.49 & 79.8    & 5.36 & 88.2\\
3 & 308 & 1.95 & 3.73    & 745  & 3.87    & 753\\\hline
\end{tabular}\quad
\begin{tabular}{ccc}
\hline
& GC & LS\\ \hline
$\eta_w$ & 99.7\% & 97.9\% \\ 
$\eta_v$ & 98.9\% & 99.7\%  \\ \hline
\end{tabular}
\end{table}

\begin{table}[!t]
    \centering

    \caption{Flow around a sphere (Exp. 4): number of nonlinear iterations $N_N$, number of linear iterations per nonlinear iteration $\overline{N}_L$, and wall-clock time $T$ in seconds. Performance metrics for the cases with $l=2$: $\eta_w$ parallel workload efficiency and $\eta_v$ vertical communication efficiency.}\label{tab:fpc}

\begin{tabular}{c|c|cc|cc}
\hline
 & \multicolumn{5}{c}{$Q_1Q_1$} \\ \hline
 & &\multicolumn{2}{c|}{GC} & \multicolumn{2}{c}{LS} \\ \hline
$l$ & $N_N$ & $\overline{N}_L$ & $T$ & $\overline{N}_L$ & $T$ \\ \hline
1 & 5  &  14.6  &  7.36   &  15.2  &  7.14\\
2 & 4  &  23.25  &  37.4   &  23  &  43.7\\ \hline
\end{tabular}\quad
\begin{tabular}{c|c|cc|cc}
\hline
 & \multicolumn{5}{c}{$Q_2Q_2$} \\ \hline
 & &\multicolumn{2}{c|}{GC} & \multicolumn{2}{c}{LS} \\ \hline
$l$ & $N_N$ & $\overline{N}_L$ & $T$ & $\overline{N}_L$ & $T$ \\ \hline
1 & 4  &  14.75  &  28.2 &  14.75  &  34.1\\
2 & 2  &  23  &  162 &  23.5  &  206\\ \hline
\end{tabular}\quad
\begin{tabular}{ccc}
\hline
& GC & LS\\ \hline
$\eta_w$ & 99.9\% & 53.8\% \\ 
$\eta_v$ & 87.1\% & 99.9\%  \\ \hline
\end{tabular}
\end{table}

In summary, both geometric multigrid algorithms require a similar number of iterations, with a slight advantage of GC over LS. It seems that the pure Dirichlet BC during local smoothing does not have a negative influence on the number of linear iterations of the solver (potentially due to the stabilization). However, a more detailed investigation is required in this direction; in particular, we plan to investigate the modified boundary conditions known in the context of domain decomposition~\cite{cai2016hybrid} or model refinement methods~\cite{tominec2025well}.

\section{Conclusions \& Outlook}\label{sec:outlook}

This work solved various scenarios of the Stokes and Navier-Stokes equations with the multigrid infrastructure of \dealii{} in parallel. We demonstrated its flexibility and modularity, which allow users to easily swap smoothers and coarse-grid solvers, choose different coarsening strategies, and add custom implementations of transfer classes. We also presented its applicability to a large set of continuous finite elements ($\bm{Q}_p/Q_{p-1}$, $\bm{Q}_{p}^d/P_{p-1}^{\mathrm{disc}}$, $\bm{Q}_p/Q_{p}$) and to $h$- and $p$-refined meshes.


The multigrid infrastructure in \dealii{} is actively maintained and developed. Current efforts focus on porting the multigrid infrastructure to GPUs, enabling us to run the simulations presented above on accelerators as well.

\section*{Acknowledgments}

The authors acknowledge collaboration with Wolfgang Bangerth, Markus Bause, Bruno Blais, and Timo Heister as well as the \dealii{} community, without which this work would not be possible.

\textbf{M.F.}\@ was partially supported by the National Science Foundation under award OAC-1835673 as part of the Cyberinfrastructure for Sustained Scientific Innovation (CSSI) program. Additional support was provided by the ERC-CZ grant LL2105 CONTACT and the P~JAC grant CZ.02.01.01/00/22\_008/0004591 FerrMion, both funded by the Czech Ministry of Education, Youth and Sports, as well as by the Czech Science Foundation project No.\@ 26-21877S and the Charles University Research Centre Program No.\@ UNCE/24/SCI/005.
\textbf{M.K.}\@ was partially supported by the German Federal Ministry of
Research, Technology and Space (BMFTR) through project “PDExa: Optimized software methods for solving partial diﬀerential equations on exascale supercomputers” (grant agreement no.~16ME0637K) and the NextGenerationEU program of the European Union.
\textbf{N.M.}\@ acknowledges funding from the European Regional Development Fund (grant FEM Poer
II, ZS?2024/06/18815) under the European Union's Horizon Europe Research and
Innovation Program and computational resources (HPC cluster HSUper) provided by the project
hpc.bw, funded by dtec.bw - Digitalization and Technology Research Center of the
Bundeswehr. dtec.bw is funded by the NextGenerationEU program of the European Union.



\section*{Code Availability}


\textbf{Exp.~1} was realized with \url{https://github.com/marcfehling/hpbox},\\ \textbf{Exp.~2} with \url{https://github.com/nlsmrg/dealii-stfem}, and \\\textbf{Exp.~3} and \textbf{Exp.~4} with \url{https://github.com/chaos-polymtl/lethe}.


\bibliographystyle{spmpsci}
\bibliography{bibliography}
\end{document}